\documentclass[11pt]{article}
\pdfoutput=1
\usepackage{mypackages}
\usepackage{mymacros}
\usepackage[framed,numbered,autolinebreaks,useliterate]{mcode}
\usepackage[affil-it]{authblk} 
\title{\LARGE Introductory Topological Data Analysis} 
\author{\large Dayten Sheffar}            
\affil{\footnotesize Department of Mathematics and Statistics\\ University of Victoria \\
Victoria, BC, V8W 2Y2, Canada\\ 
dsheffar@uvic.ca}  
\date{\vspace{-5ex}}
\begin{document} 
\maketitle 
\section{Abstract}
This paper introduces topological data analysis. Starting from notions of a metric space and some elementary graph theory, we take example sets of data and demonstrate some of their topological properties. We discuss simplicial complexes and how they relate to something called the Nerve Theorem. For this we introduce  notions from the field of topology such as open covers, homeomorphism and homotopy equivalences. This leads us into discussing filtering data and deriving topologically invariant simplicial complexes from the underlying data set. There is then a small introduction to persistent homology and Betti numbers as these are useful analytical tools for TDA. An accompanying online appendix for the code producing the bulk of the figures in this paper is available at bit.ly/TDA\_2020.
\begin{multicols}{2}
\section{Introduction}
Data is an increasingly important resource in our world. The complexity of data that can be worked with grows as innovative analytical techniques proliferate in private industry and academia. 

One such method that has been developed is called Topological Data Analysis (TDA), a novel tool that allows investigating complex and high-dimensional data by forming lower dimensional approximations that preserve structure and connectivity-information in the data \cite{carlsson}. Suppose we perform an experiment and collect some data. We want to look at it and infer relationships the features might have, whether they possess a linear relationship, an inverse relationship, and so forth. 

\begin{figure}[H]
    \centering
    \begin{minipage}{0.495\linewidth}
        \includegraphics[width=0.95\linewidth]{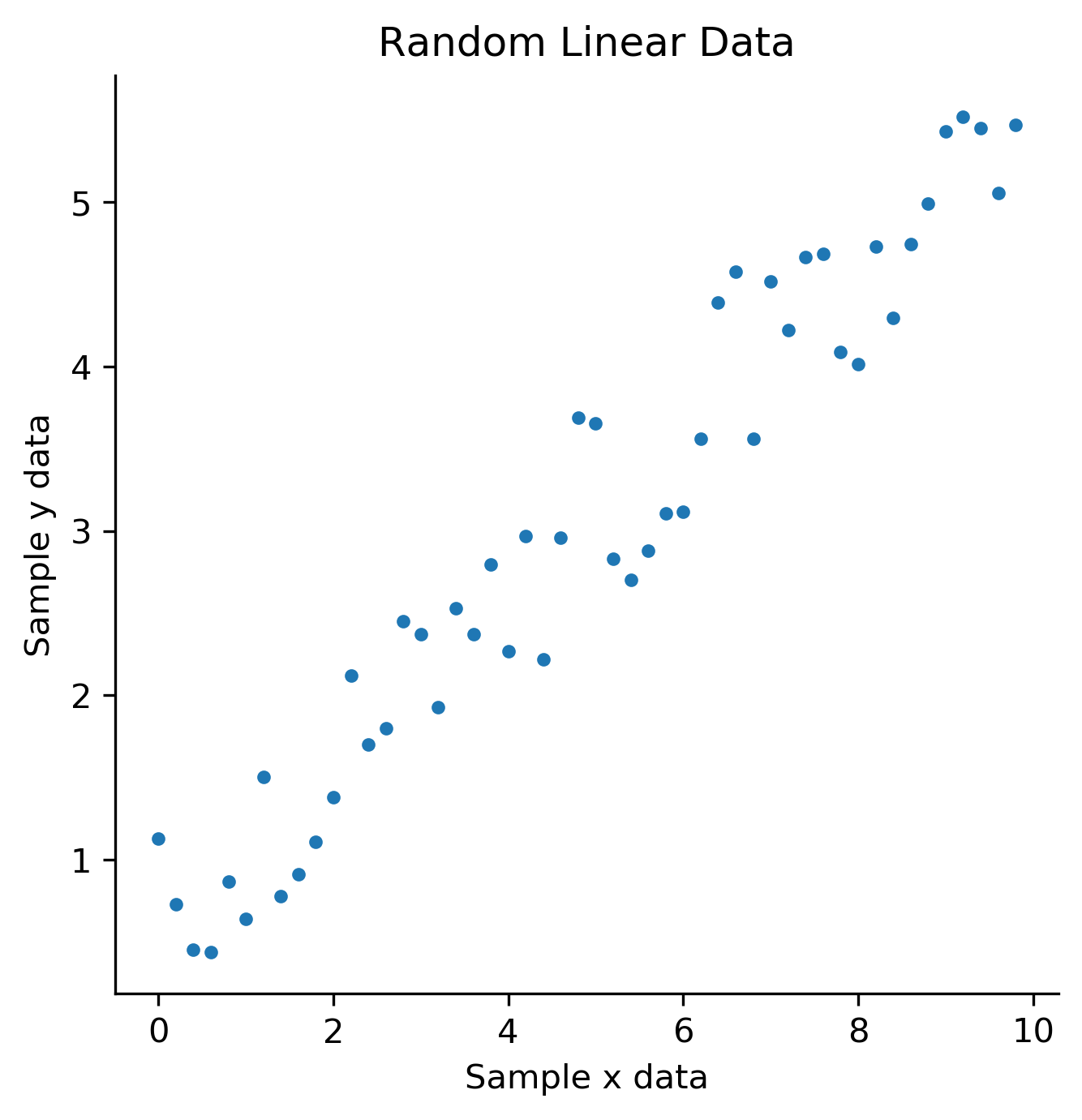}
        \caption{}
        \label{lin_plot}
    \end{minipage}
    \begin{minipage}{0.495\linewidth}
         \centering
        \includegraphics[width=1.00\linewidth]{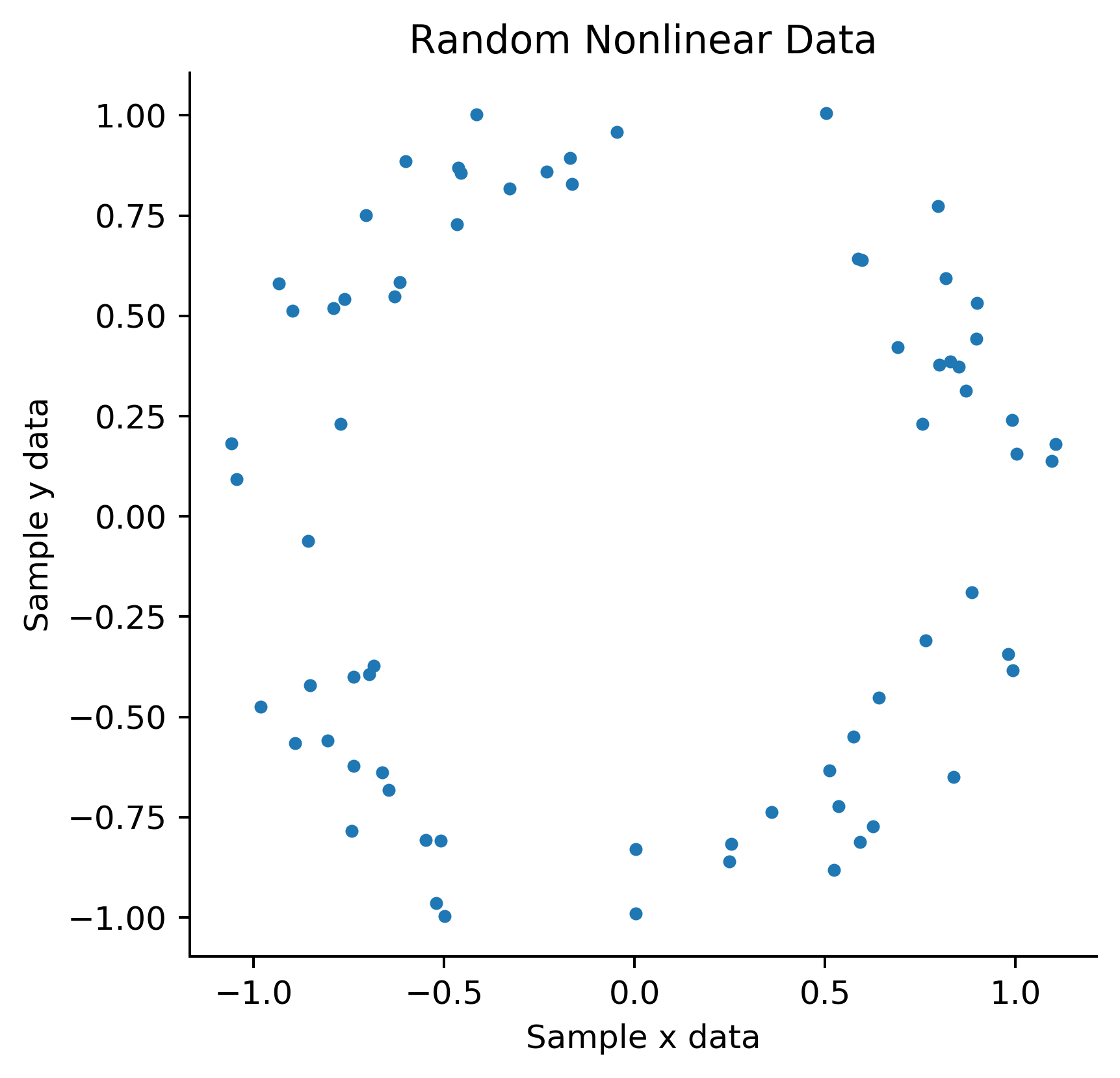}
        \caption{}
        \label{circ_plot}
    \end{minipage}
    \caption*{\textit{Left we have some randomly spread linear data, and right some random annular data  } $\mathcal{A}$}
\end{figure}

Perhaps the linear data depicted in figure \ref{lin_plot} could come from a simple experiment in a lab. One could then perform say, a linear regression on the data (the line of best fit) to find the relationship between the features of the experiment. 

Linear data such as this have a fairly simple distribution, we turn to the somewhat more interesting case of figure \ref{circ_plot}, that of annular data we will call $\mathcal{A}$, and one will be able to quickly infer the results of this paper back to this case. 

Perhaps on $\mathcal{A}$ the idea would be to do a non-linear fit of a circle to the data by minimizing the euclidean distance between the curve and the data. While this can be a useful exercise, ideas have emerged from topology for gaining novel insights into your data regarding its geometric shapes, connectivity, and its similarity to another data set.

    \section{Topology \& Graph Theory}

Very quickly, a \textbf{topology} on our a $\mathbb{X}$ is a collection of subsets of $\mathbb{X}$ that we call $\mathcal{T}$. We need $i)$ that the empty set (having nothing) and $\mathbb{X}$ itself (having everything) are both in the topology $\mathcal{T}$, $ii)$ the union of the elements of any subcollection of $\mathcal{T}$ is in $\mathcal{T}$, and $iii)$ the intersection of the elements of any finite subcollection of $\mathcal{T}$ is in $\mathcal{T}$.

We now need to establish many fundamental notions before we can get into the bulk of what topological data analysis is.

\subsection{Metric Topology}
We will be considering the Euclidean plane since our data in $\mathcal{A}$ is a subset of $\R^2$, which has something called a metric topology. If you're familiar with \textit{Euclidean distance}, that's the notion we wish to generalize here.

We can define a \textbf{distance function} $d$ that takes arguments $(x,y)$ from $\R \times \R$ to $\R$ which is $i)$ always non-negative unless $ii)$ the distance is zero which is if and only if we're looking at the distance from a point to itself $d(x,x)=0$. We need  $iii)$ that the distance is symmetric, denoted $d(x,y) = d(y,x)$ when $x \not= y$, and finally $iv)$ it obeys the triangle inequality 
$$ d(x,y) + d(y,z) \geq d(x,z) $$ for any $x,y,z$ in $\R$. We often denote the metric topology by $(\mathbb{X},\mathcal{T}_{\vert\cdot\vert})$

\subsection{$p$-Norms}
A $p$-norm is a notion of the length of a vector or the distance from the origin to a point in a metric topology \cite{munkres}, and for $x \in \R^n$ it is denoted by
$$ L_p = \norm{\mathbf{x}} = \left(\sum_{i=1}^n (x_i)^p \right)^{-p} $$
and if we're in $\R^2$ this simply reduces to the common Pythagorean notion of distance denoted by $L_2$, or the classical $c = \sqrt{a^2 + b^2}$.

If we take $p=2$ and $x_n, y_n \in \R^n$, we can state the Euclidean distance between two vectors or points \cite{munkres}, 
$$ d(x,y) = \norm{\mathbf{x}-\mathbf{y}} = \left(\sum_{i=1}^n  (x_i - y_i)^2 \right)^{-2}$$
and later we will need a derivative notion of the \textit{infinity norm}, denoted $L_\infty$, which is the maximum norm in the point cloud. This would merely indicate the point furthest from the origin, and what we need however is the maximum distance between any two points of $\mathcal{A}$ for which we will we borrow the label $L_\infty$ for brevity,

$$ L_\infty = \max_{x_i,x_j\in\mathcal{A}} \left\{ \norm{\mathbf{x}_i-\mathbf{x}_j}\right\}. $$ 

Of note (but not proven here) are the facts that the points constituting this maximum distance are part of the convex hull of the data (see figure \ref{fig:convex_hull}. Imagine putting a noose outside the boundary of your data set and tightening it until you have straight lines between the most exterior points. The pair constituting the maximum are found with a method called rotating calipers (not investigated here), once the convex hull is found. Rotating calipers is a technique from computational geometry that compares the distance between all points of the convex hull and takes the maximum value found as the $L_\infty$ norm.

\subsection{Abstract Simplexes and the Convex Hull}

We need to distinguish a simplex from simplicial complexes, which will be covered later in section \ref{simplicial_complexes}. Given a set (we can use $\mathcal{A}$), a \textbf{simplex} is the complete graph induced on the set (this can happen in $n$-dimensions). Just think of connecting every single point of $\mathcal{A}$ with an edge. If the size of $\vert \mathcal{A} \vert = k+1$ then we obtain a $k$-dimensional abstract simplex. The $k+1$ points are predictably called the vertices of the simplex and the \textbf{simplices} are the faces obtained by joining these vertices.

\begin{figure}[H]
    \centering
        \includegraphics[width = 0.85\linewidth]{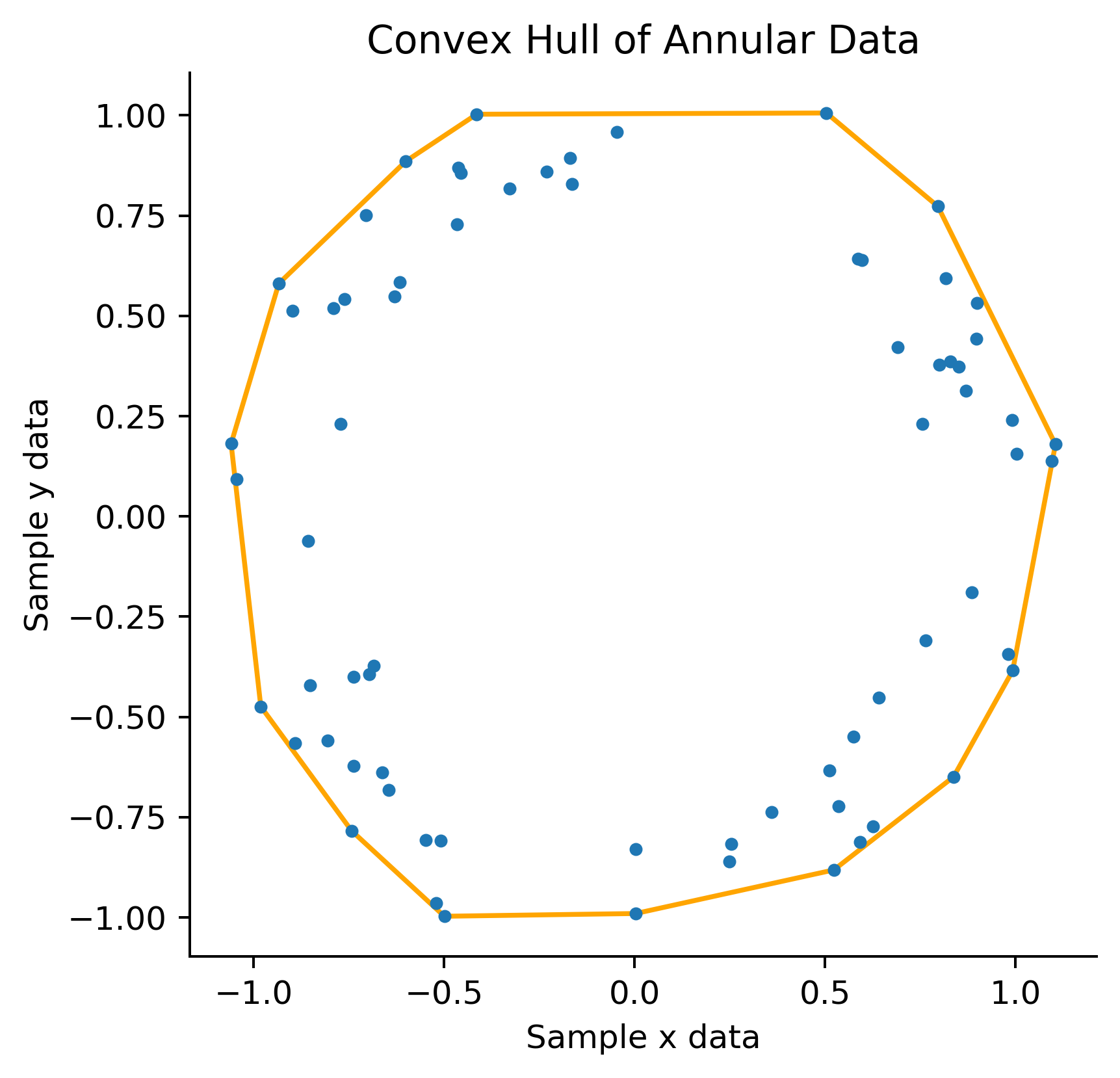}
        \caption{\textit{Convex Hull of $\mathcal{A}$}}
        \label{fig:convex_hull}
\end{figure}

For our case, the \textbf{convex hull} (see figure \ref{fig:convex_hull}) is the minimal convex subset of $\R^2$ consisting of vertices and straight edges that contains $\mathcal{A}$. Finding the convex hull is useful in itself as it is the boundary of your data values, which can provide guidance in further testing/data collection, or aid in decision optimization such as seeking to produce data in only certain sub-regions.

\subsection{Open Ball}
We also want to introduce the important notion of an \textbf{open ball} of radius $r$ around some point $x_0$ in $\mathcal{A}$, denoted as
$$ \mathcal{B}_d(x_0, r) = \{ y \vert d(x,y) < r \} $$
as being the distance within a radius of a point not including the boundary \cite{munkres} (we will relax this boundary condition later). 

Hence, a topology is induced on $\mathcal{A}$ by our distance function $d$. We have the metric topology on $\R$ denoted ($\R,\mathcal{T}_{\vert\cdot\vert})$, and examine the subspace topology of balls centered at points of $\mathcal{A}$. A subspace topology is simply a set $\mathbb{X} \subset \R^2$ together with the topology $\mathcal{T}_2$ found by taking intersections with open sets of the metric topology, specifically, $\mathcal{T}_2 = \{ \mathbb{X} \cap U \vert U \in \mathcal{T}_{\vert\cdot\vert} \} $ \cite{munkres}.

Let us now turn to an examination of how the data from $\mathcal{A}$ connects with itself as we expand the radius of these balls centered at points of $\mathcal{A}$.
    \section{Connectivity of the data} \label{connectivity}
The Euclidean plane with the metric topology is said to be connected. \textbf{Connectedness} of a topological space indicates there does not exist what is called a separation on the space. A \textbf{separation} is when you can find non-empty, disjoint open subsets $U$ and $V$ of some set $\mathbb{X}$ imbued with a topology. These will come up shortly.

\begin{figure}[H]
    \centering
    \begin{minipage}{0.48\linewidth}
        \includegraphics[width = 1\linewidth]{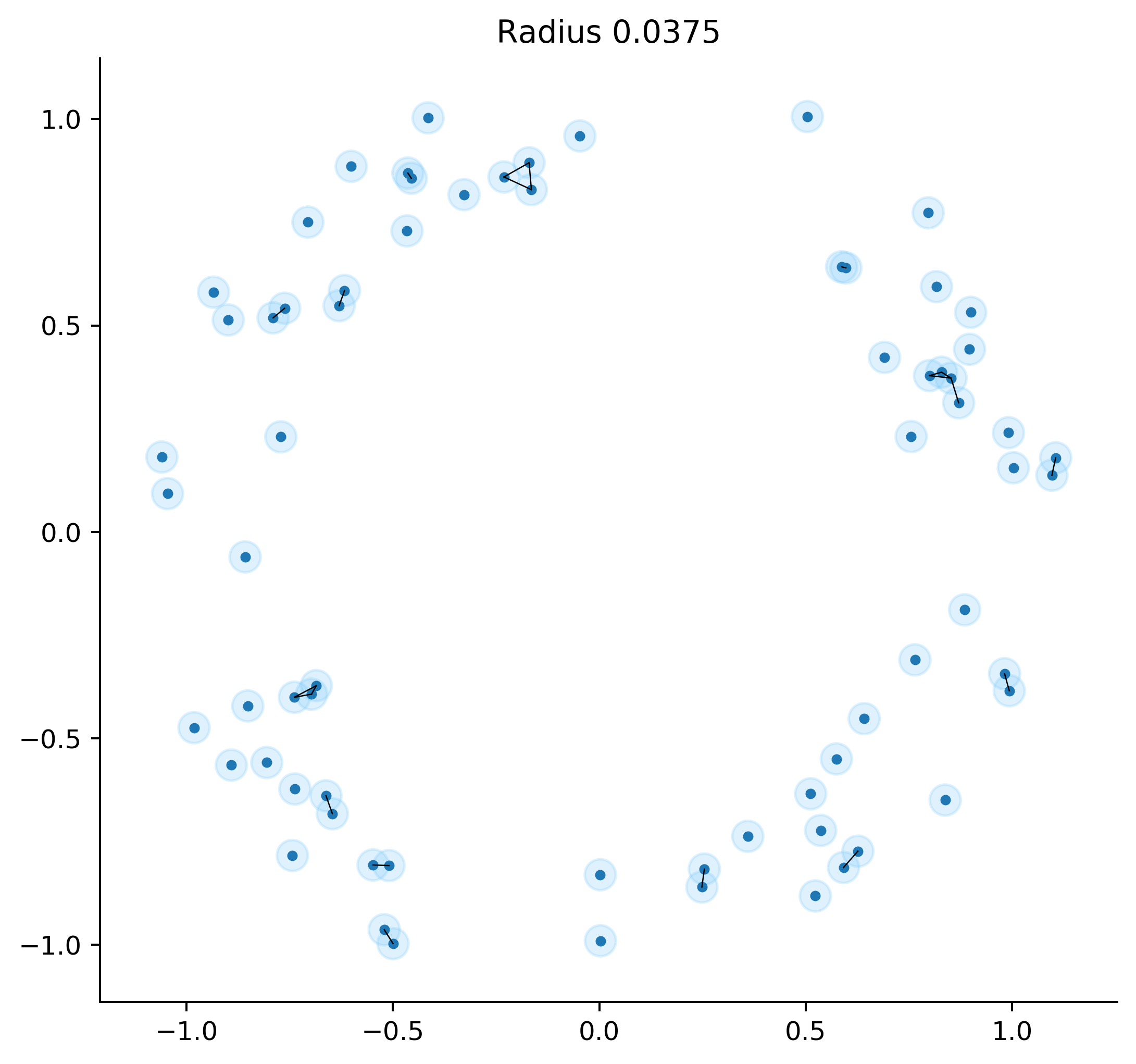}
        \caption{}
        \label{fig:rad0}
    \end{minipage}
    \begin{minipage}{0.48\linewidth}
        \centering
        \includegraphics[width = 1\linewidth]{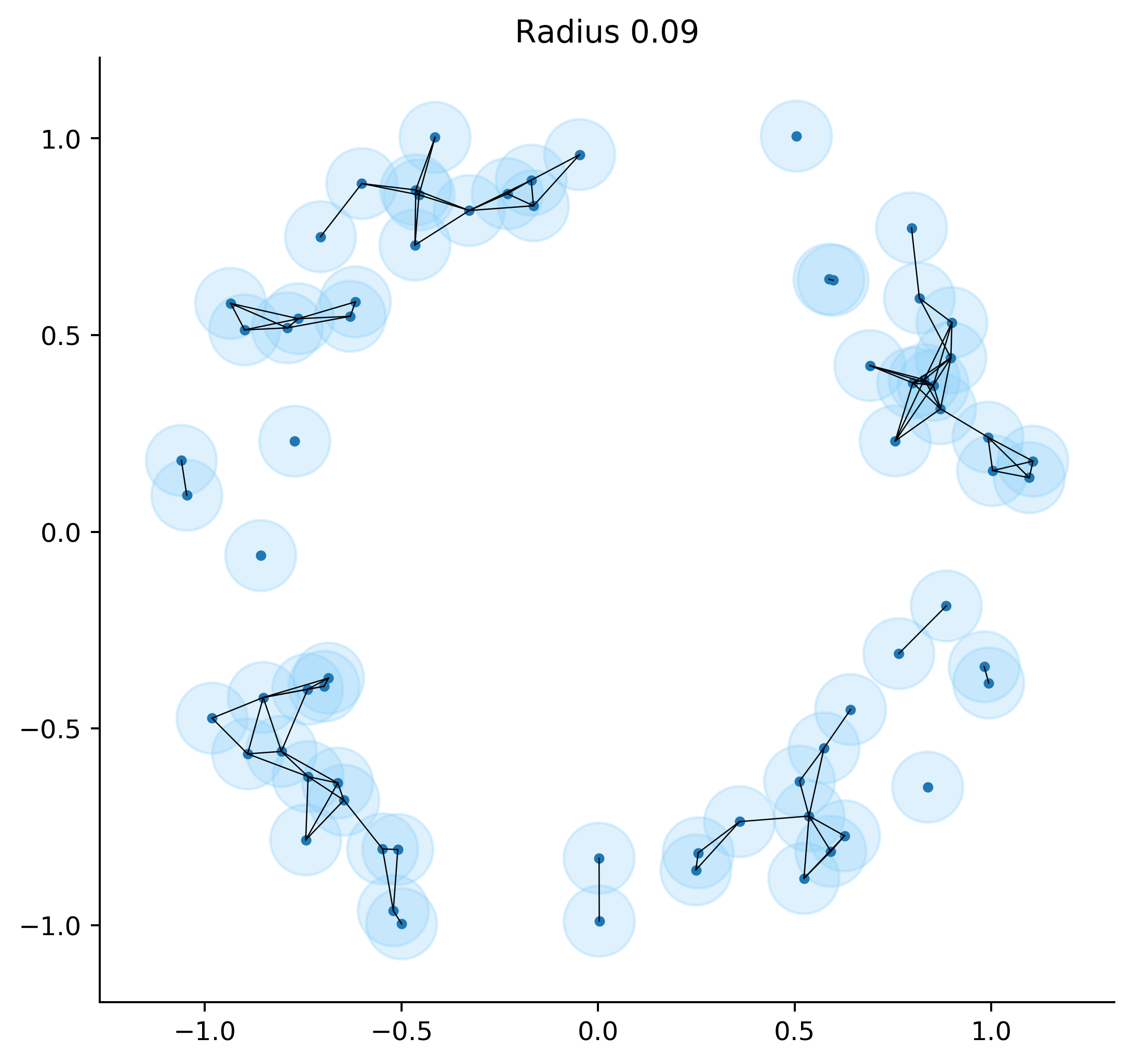}
        \caption{}
        \label{fig:rad1}
    \end{minipage}
    \caption*{\textit{Initial radius on the left, increasing connections hint at denser clusters on the right}}
\end{figure}

Let's take $\mathcal{A}$ and around each point $x_i$ draw a ball $\mathcal{B}(x_i, r)$ for some small initial radius $r$ (figure \ref{fig:rad0}). We want to take snapshots for various increasing $r$ and get an idea of how $\mathcal{A}$ is connecting for various radii.

What we will do is make a note whenever the closures of $\mathcal{B}(x_i,r)$ and $\mathcal{B}(x_j,r)$ intersect for the first time on a point in the plane. This is where we could relax our open ball constraint, for example, the closure of the open ball around $x_i$ is merely the closed ball around $x_i$ of radius $r$. We say that $x_i$ and $x_j$ are $2r$ apart in distance, connect them with an edge, and visualize our graph \ref{fig:rad0}.

\begin{figure}[H]
    \centering
        \includegraphics[width = 0.85\linewidth]{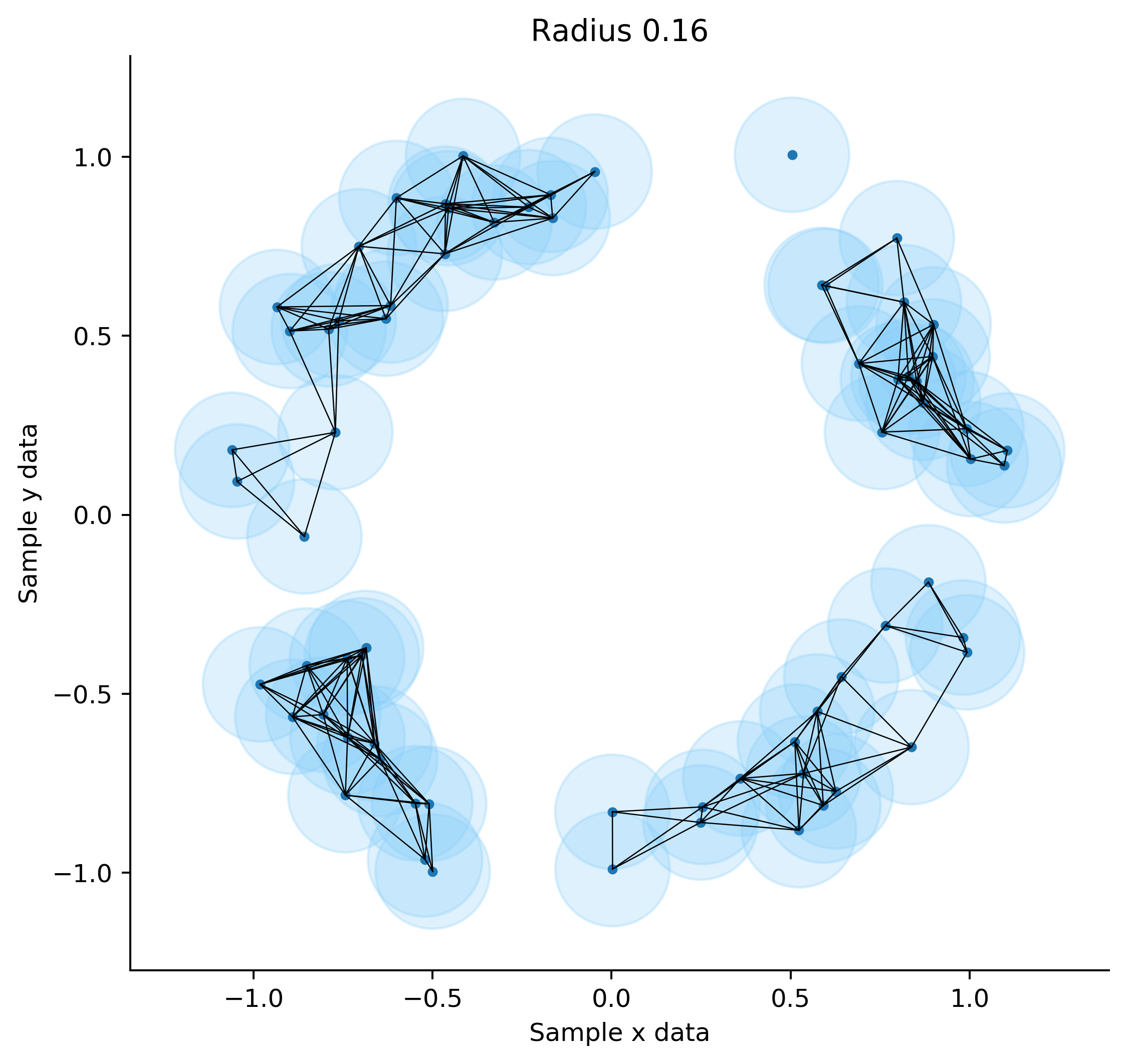}
        \caption{\textit{Clusters in $\mathcal{A}$}}
        \label{fig:rad2}
\end{figure}

First, we have sparse groupings of vertices connected by edges. As we increase $r$, more edges begin to appear in small connected components in the graph that we call clusters (figure \ref{fig:rad2}), the resulting figure may also be called a nearest-neighbor-graph (for a given radius). Herein, we begin to note the first signs of a structure on the data, sometimes called it's connectivity-information \cite{carlsson}. Typically, clusters denote important categorical information about the underlying features/variables of the experiment (perhaps when some initial conditions of your experiment are met, the data points are more likely to be found clumped in certain ways). 

Next, significant clusters may start to emerge that vary from clusters achieved with smaller radii, typically subsuming them. This raises an important consideration for what level of granularity or resolution you want in your analysis, or how many clusters would you like?. In topology we consider this as the \textbf{fineness} of a topology \cite{munkres}. Note that for the clusters themselves, a separation can be found, indicating that the we have not yet found the radius for which $\mathcal{A}$ is wholly connected.

As the radius of the balls expand, we eventually obtain a single component for the first time, and we can say that $\mathcal{A}$ is finally a single connected component. The figure can be a chain, have loops or holes, be a loop, have multiple flares, or any other structures. We will discuss this later in section \ref{betti_numbers} with what are called Betti numbers. 

\begin{figure}[H]
    \centering
    \begin{minipage}{0.48\linewidth}
        \centering
        \includegraphics[width = \linewidth]{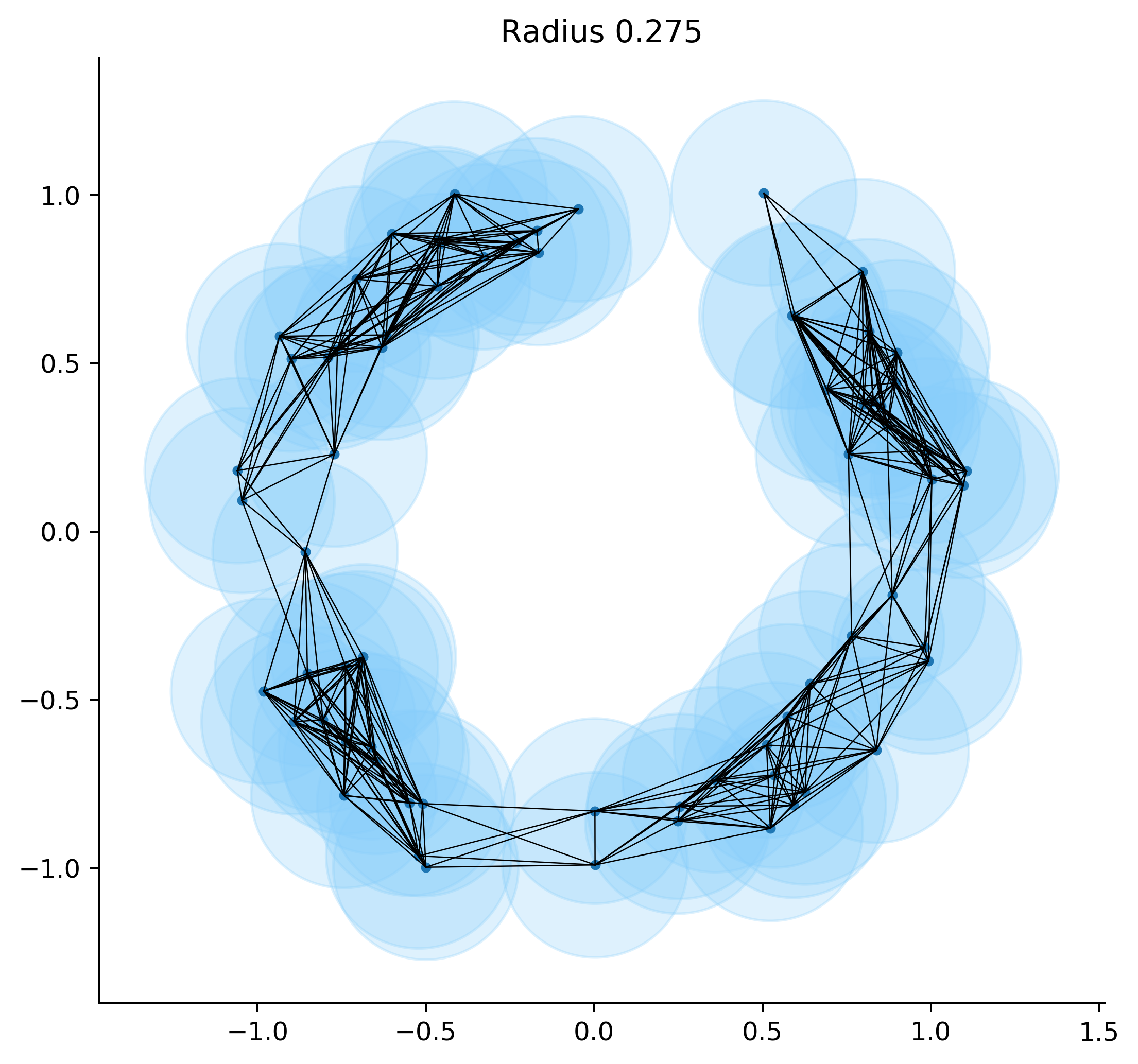}
        \caption{}
        \label{fig:rad3}
    \end{minipage}
    \begin{minipage}{0.48\linewidth}
        \centering
        \includegraphics[width = \linewidth]{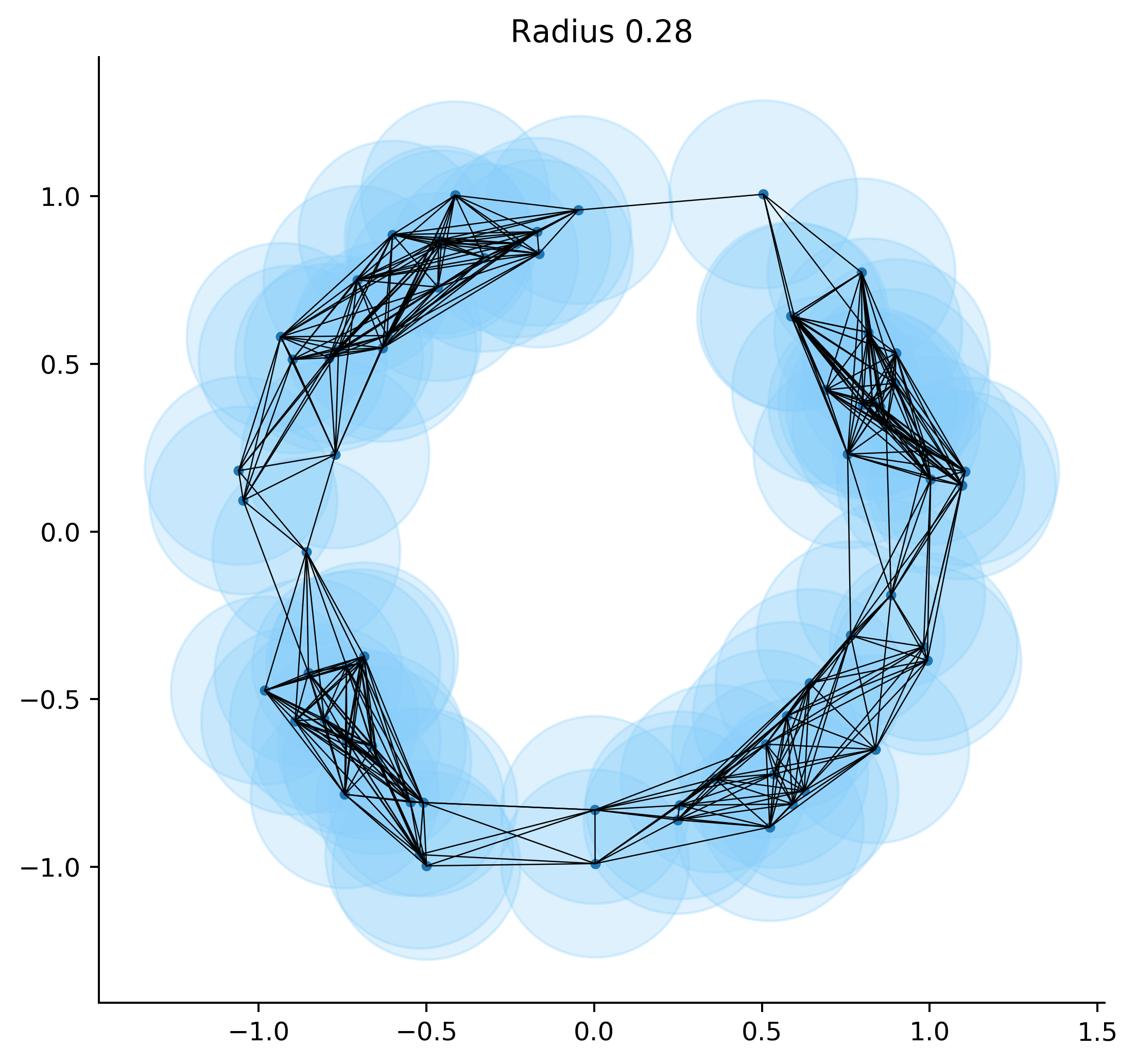}
        \caption{}
        \label{fig:rad4}
    \end{minipage}
    \caption*{\textit{Largest singular component of $\mathcal{A}$ on the left, before increasing $r$ and first obtaining the loop structure on the right for some $r^*$ such that $0.275 < r^* < 0.28$}}
\end{figure}

Following this, still increasing the radius $r$, we begin to see the loop in the data emerge as in figure \ref{fig:rad4}. For data that is not randomly generated, there could be some substantial reason as to the existence of a cycle or holes.

So as $r$ approaches $r^*$ the balls around the respective $x_i$ and $x_j$ intersect and they obtain the final edge to complete the loop for the first time. Once we've reached this radius, we see that increasing the radius yields no further structure for quite some time. This range begins for a radius $r > r^*$ as seen between figures \ref{fig:rad3} and \ref{fig:rad4}. We might safely say we've obtained the characteristic descriptor of this data (in this case being annular/having a loop) and may be comfortable stopping our analysis at this stage. The Betti numbers seen later with persistence diagrams characterize what radii reveal what structures (see section \ref{betti_numbers}).

\begin{figure}[H]
    \centering
    \begin{minipage}{0.48\linewidth}
        \centering
        \includegraphics[width = \linewidth]{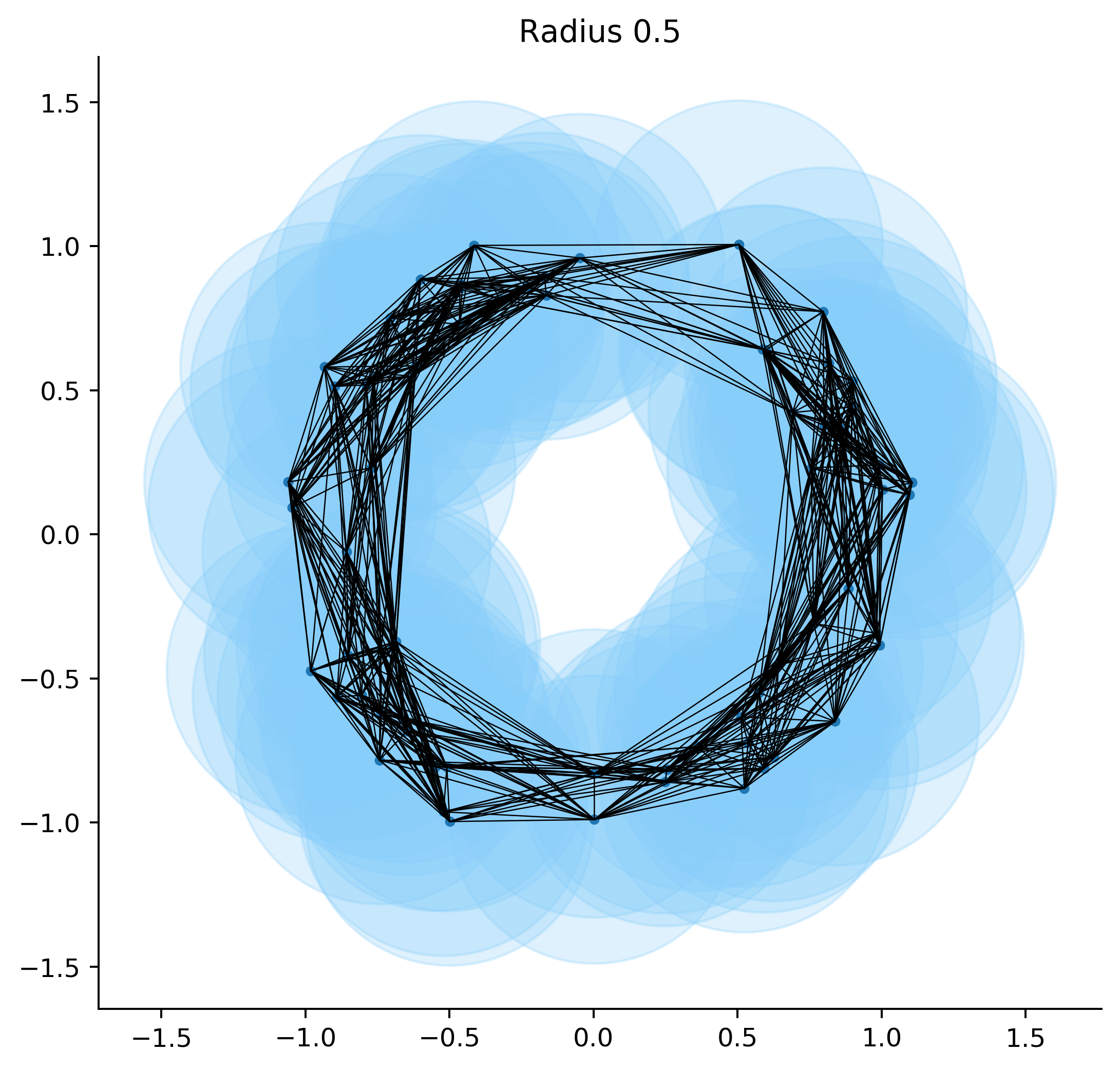}\caption{}
        \label{fig:rad5}
    \end{minipage}
    \begin{minipage}{0.48\linewidth}
        \centering
        \includegraphics[width = \linewidth]{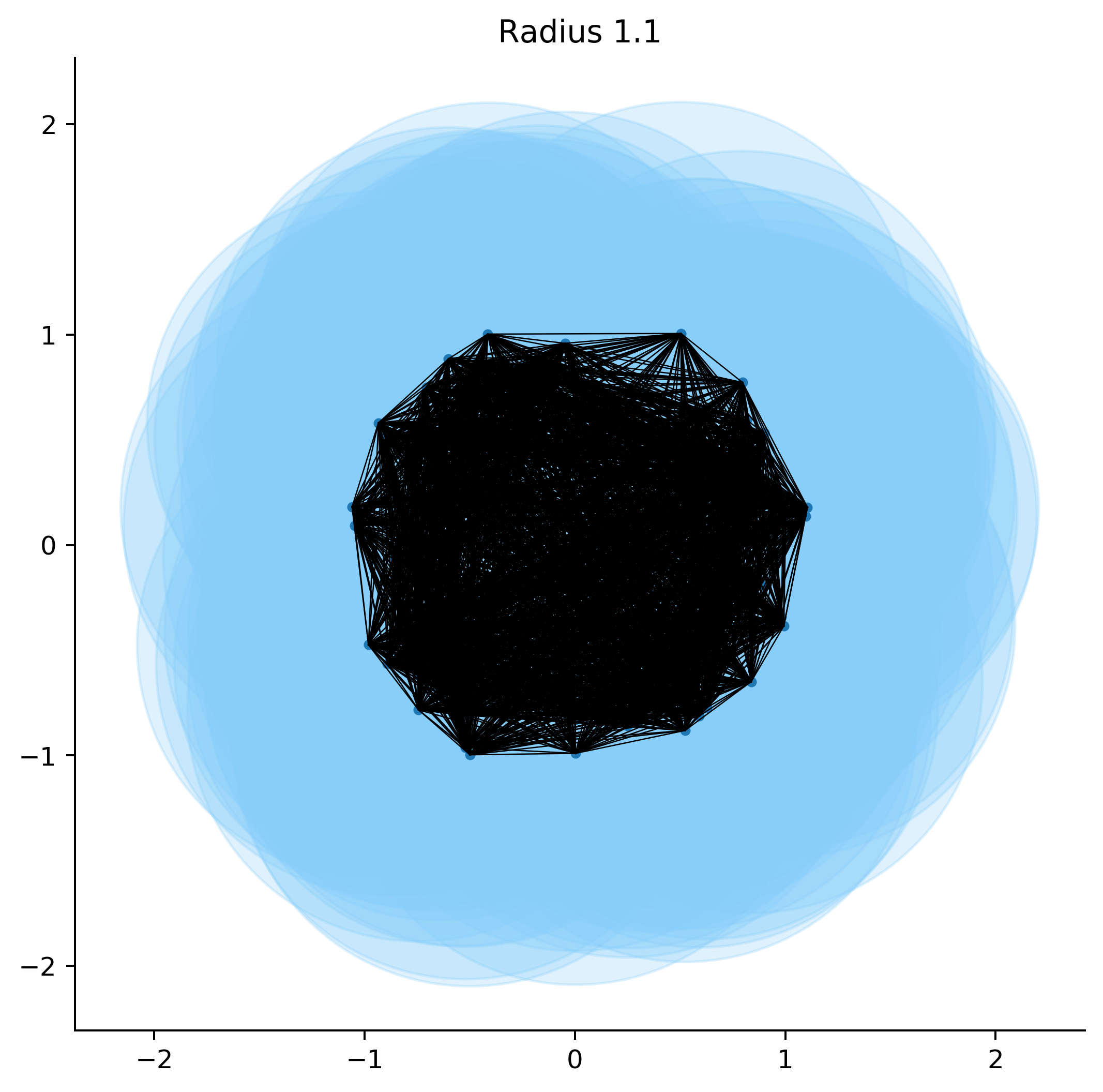}
        \caption{}
        \label{fig:rad6}
    \end{minipage}
    \caption*{\textit{Left, there is no new macro-structure revealed for $r > r^*$ for some time. Right, the next development then seems to be in obtaining the convex hull}}
\end{figure}

However, as mentioned before, for a ball radius extended past the $L_\infty$ norm as described, we collect the convex hull for free. In our example, it is approximately circular, reflecting some noise in the random genesis of our data.

\subsection{Simplicial Complexes} \label{simplicial_complexes}

With a better understanding of the process of connecting our data, we return to the topic of simplicial complexes, where there are two standard ways we will see of interpreting the constructions in section \ref{connectivity} where we expanded balls around the points of $\mathcal{A}$.

First, let's discuss the \textbf{Vietoris-Rips complex}, an abstract simplicial complex which does not necessarily permit a geometric realization in the euclidean plane (which is precisely the case for our data $\mathcal{A}$) \cite{chazal}. This was discussed earlier as the $k$-dimensional simplex in $\R^2$ for our case.

\begin{figure}[H]
    \centering
        \includegraphics[width = 0.85\linewidth]{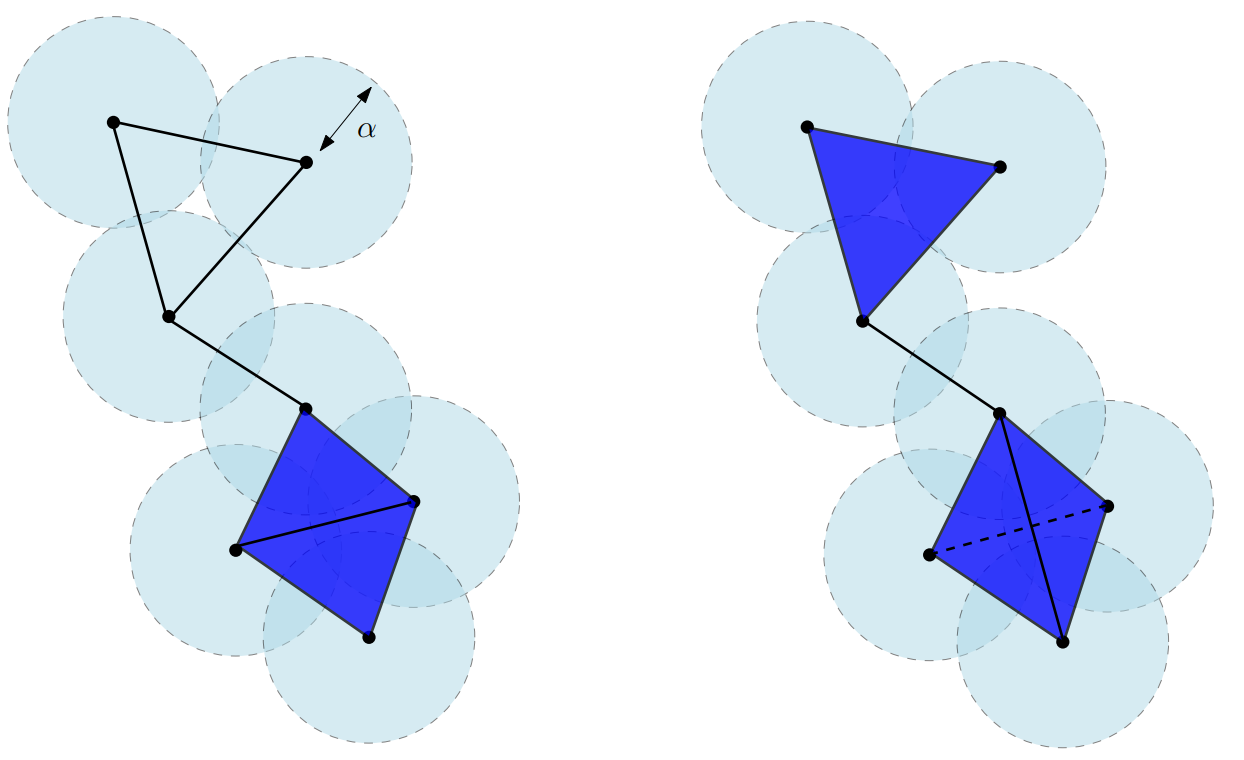}
        \caption{\textit{\v{C}ech (left) and Rips (right) complexes on a toy dataset $\mathbb{X}$} - figure via Chazal and Michel (2017)}
        \label{fig:simplicial_complexes}
\end{figure}

 Figure \ref{fig:simplicial_complexes} elucidates the fact that the Rips complex can be 3-dimensional in a 2-dimensional space (and is thus not embedded in the plane). Denoted Rips$_\alpha(\mathcal{A})$ for a fixed radius $\alpha$, then $d(x_i,x_j) \leq \alpha$ for all (i,j) \cite{chazal}. 
 
 The \textbf{\v{C}ech complex}, \v{C}ech$_\alpha(\mathbb{X})$ is the set of simplexes on some $k+1$ vertices, such that the balls $\mathcal{B}(x_i,\alpha)$ have non-empty intersection \cite{chazal}. The \v{C}ech complex has dimension 2 as in figure \ref{fig:simplicial_complexes}, and is thus embedded in $\R^2$,

We note two facts, first, the methods are related by
$$ \text{Rips}_\alpha(\mathbb{X}) \subset \text{\v{C}ech}_\alpha(\mathbb{X} \subset \text{Rips}_{2\alpha}(\mathbb{X})  $$
and further, both methods are based on the same vertices and edges \cite{chazal} . They differ in how they treat the collection of faces. The simplices in figure \ref{fig:simplicial_complexes} (colored in blue) for the \v{C}ech$_\alpha$ complex are the union of two adjacent triangles, but for the Rips$_{2\alpha}$ complex we obtain a 3-dimensional tetrahedron \cite{chazal}.

The \v{C}ech complex may be computationally expensive, and so, among other algorithmically efficient variants that exist, the Rips complex is often used for more expensive to compute data sets \cite{bubenik}.

\subsection{Homeomorphic and Homotopy Equivalence}

We need to press the brakes for a minute and think about where we might like to go from here. We'd like a notion of how to check if two topological spaces $\mathbb{X}$ and $\mathbb{Y}$ share similar topological features. One method of proceeding would be to check if $\mathbb{X}$ and $\mathbb{Y}$ are \textbf{homeomorphic}, that is, if there exists bijective functions $f : \mathbb{X} \to \mathbb{Y}$ and $g: \mathbb{Y} \to \mathbb{X}$ such that $f\circ g$ and $g \circ f$ are the identity maps for $\mathbb{Y}$ and $\mathbb{X}$ respectively \cite{chazal}. Then $f=g^{-1}$ and $f$ is said to be a \textbf{homeomorphism} \cite{munkres}. 

However, this may be too strong a notion \cite{chazal} to ensure the similarity we seek between the spaces. Why we're doing this will become apparent in the next few sections. First, if homeomorphic is too strong, relaxing the conditions leads to the related notion of being homotopic. 

Two continuous maps $f_0$ and $f_1$ from $\mathbb{X}$ to $\mathbb{Y}$ are called \textbf{homotopic} if there is a continuous mapping $h: \mathbb{X} \times [0,1] \to Y$ such that $h(\mathbb{X},0) = f_0(\mathbb{X})$ and $h(\mathbb{X},1) = f_1(\mathbb{X})$ \cite{munkres}. Note, we say that $\rho$ is a parameter for $h(x,\rho)$. Then $\mathbb{X}$ and $\mathbb{Y}$ are \textbf{homotopy equivalent} if there exist $f$ and $g$ such that $f\circ g$ and $g\circ f$ are homotopic to the identity map of $\mathbb{Y}$ and $\mathbb{X}$ respectively, and $f$ and $g$ are said to be homotopy equivalent as well \cite{chazal}. Homotopy equivalence is weaker than homeomorphic, as desired, however the spaces $\mathbb{X}$ and $\mathbb{Y}$ will retain many topological invariants for us to study \cite{chazal}, indeed they share the same homology (section \ref{betti_numbers}). Lastly, being homeomorphic implies being homotopy equivalent but the converse does not hold.  

As we will need this, let us briefly state that a space is \textbf{contractible} if it is homotopy equivalent to a point \cite{chazal}. Balls and convex sets in $\R^d$ are basic examples of contractible spaces.

\subsection{Covers and Index Sets}
Briefly, a collection of subsets of a space $\mathbb{X}$ is a \textbf{covering} of $\mathbb{X}$ if the union of elements of the collection equals $\mathbb{X}$. We call the collection an \textbf{open cover} if the subsets are open subsets of $\mathbb{X}$.

A set $\mathcal{I}$ is called an \textbf{index set} if its members label another set $\mathbb{X}$. It is enough to think of it as a surjective (onto) mapping from $\mathcal{I}$ to $\mathbb{X}$, and we can write the indexed collection as $\{\mathbb{X}_i\}_{i \in \mathcal{I}}$.

\section{Nerve Theorem}
Altogether, we arrive at the Nerve Theorem, of which there are several formulations. Here we provide a fairly accessible formulation of the theorem, more sophisticated versions exist.

\smallskip

\noindent \textit{\textbf{The Nerve Theorem}} Let $\mathcal{U} = (U_i)_{i\in \mathcal{I}}$ be a cover of a topological space $\mathbb{X}$ by open sets
such that the intersection of any subcollection of the $U_i$’s is either empty or contractible. Then,
$\mathbb{X}$ and the nerve $C(\mathcal{U})$ are homotopy equivalent. \cite{chazal}

\smallskip

\noindent\textbf{\textit{Proof}} Convex sets in $\R^2$ are contractible. Take a ball $\mathbb{X} = \mathcal{B}(x_0,r)$ and a point $\mathbb{Y} \in \mathbb{X}$, two spaces in $\R^2$. We show the two spaces are homotopy equivalent. First, if $\forall x \in \mathbb{X}$ we let $f(x)=\mathbb{Y}$ and $g(\mathbb{Y}) \in \mathbb{X}$, then $(f\circ g)(\mathbb{Y})=\mathbb{Y}$, and therefore $f\circ g \simeq I_\mathbb{Y}$. For $g\circ f$, we now have $(g\circ f)(x)\in \mathbb{X}$. Retracting $\mathbb{X}$ along radial lines shows  $g\circ f \simeq I_\mathbb{X}$. If we want to write down a homotopy between the two maps, then 
$$ h(x,\rho) = (1-\rho)g(\mathbb{Y}) + \rho x $$
is an example that works. That is, if $\rho=0$ we get back $g(\mathbb{Y})$, and if $\rho=1$ we get back $x$.

Consequently, granted $\mathcal{U} = (U_i)_{i\in\mathcal{I}}$ is a collection of convex subsets of $\R^d$, then $C(\mathcal{U})$ and $U_i\in\mathcal{I}$ are homotopy equivalent. Or, given $\mathbb{X}$ is a set of points in $\R^d$, then the \v{C}ech complex \v{C}ech$_\alpha$($\mathbb{X}$) is homotopy equivalent to the union of balls 
$$\bigcup_{x \in\mathbb{X}} B(x, \alpha)$$ for some fixed radius $\alpha$. \hfill $\square$

This is quite remarkable in that a much lower dimensional approximation of a set of points/data retains topological invariants such as loops and holes.
    \subsection{Filtering the data}

We want to illustrate the Nerve Theorem to gain some intuition about its components. We're going to introduce a denser set of annular data $\mathcal{A}'$, and introduce an arbitrary coloring to aid in seeing exactly what is occurring.

\begin{figure}[H]
    \centering
    \includegraphics[width=0.85\linewidth]{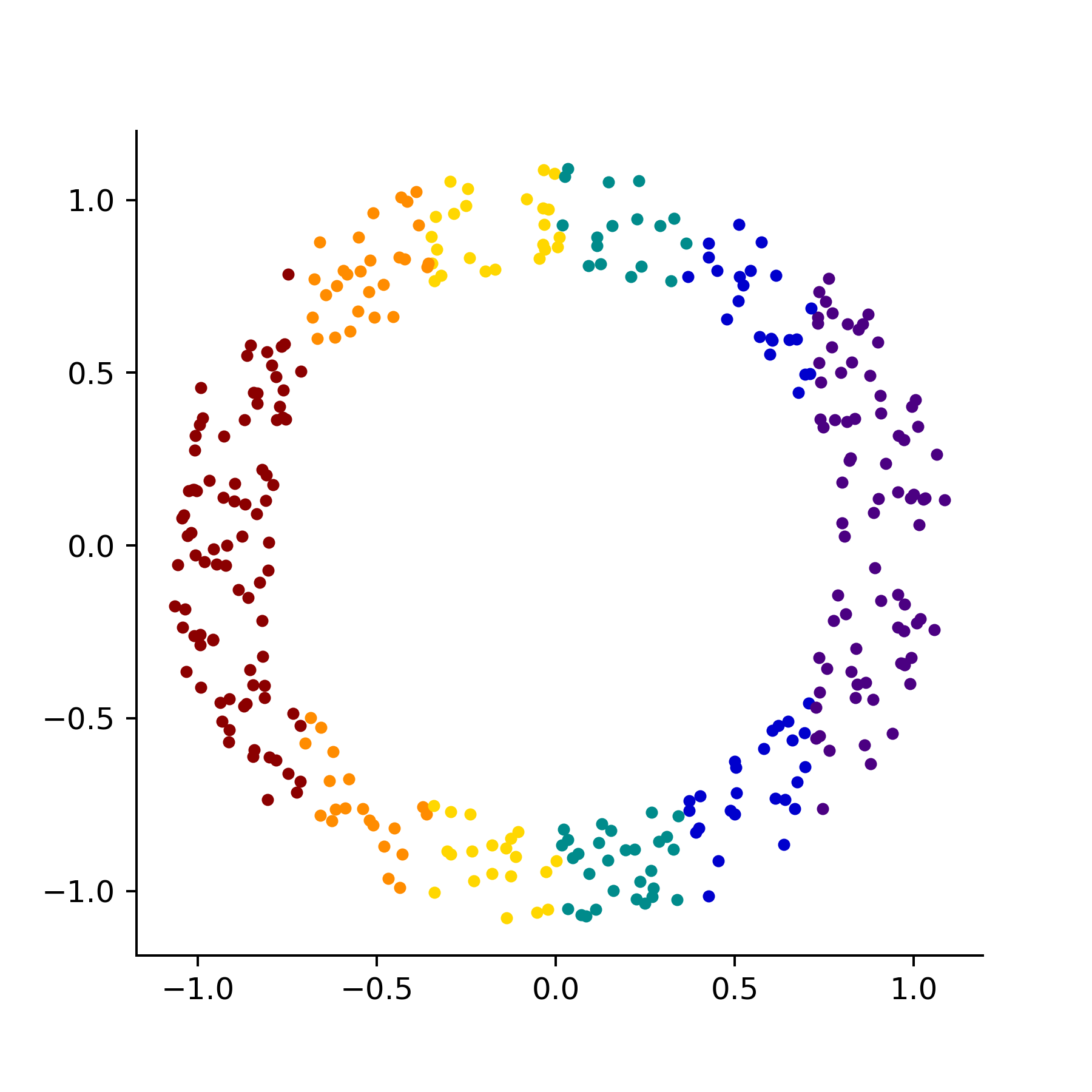}
    \caption{\textit{Partition of} $\mathcal{A}'$ \textit{(binning the data)}}
    \label{fig:filter}
\end{figure}

We want to filter or \textit{bin} the data \cite{singh}. For 2-dimensional objects we could bin along an axis, a sloped line, any clusters we have found or perhaps a curved or even non-linear (i.e. a circular) curve if we so sought. In 3-dimensions the premise is analogous, perhaps binning a 3d shape by its height. 

\begin{figure}[H]
    \centering
    \includegraphics[width=0.95\linewidth]{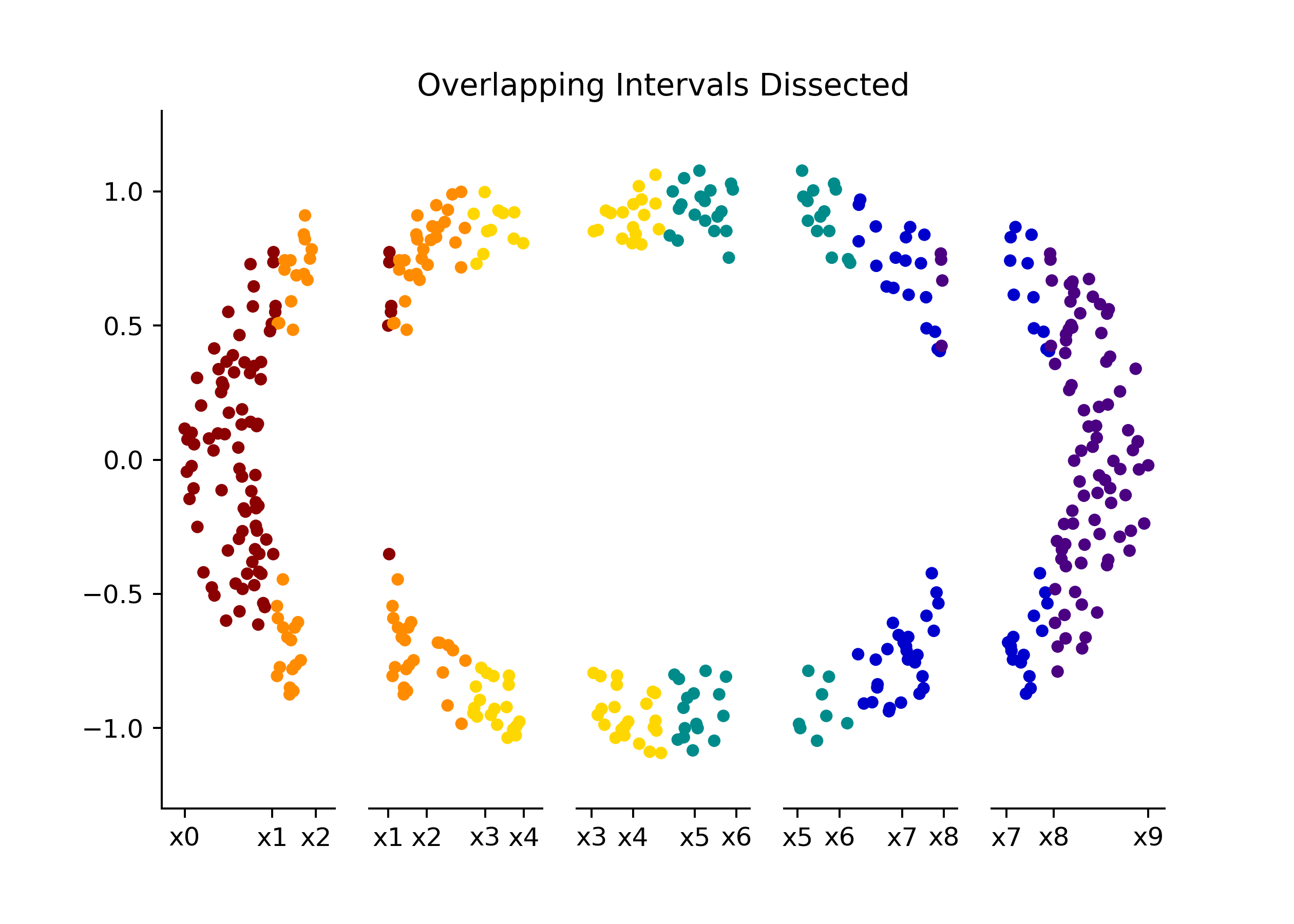}
    \caption{\textit{Overlapping intervals from } $\mathcal{A}'$}
    \label{fig:overlaps}
\end{figure}

We will bin $\mathcal{A}'$ from left to right, and what we want are overlapping intervals along the $x$-axis (so they have non-empty intersections). Conceptually it looks something like figure \ref{fig:overlaps}. Notice, the red line interval spans all the red points of $\mathcal{A}'$, but \textit{also} some orange points. Then the orange line interval spans some red points, through the orange points, and even into some yellow points.

This overlap is precisely the key factor in binning data, we're not merely taking the red points as our interval of interest, we want that overlap into the orange, and et cetera. It might seem like some double counting is occurring here, and you'd be right in thinking so! This has been made explicit in figure \ref{fig:overlaps}, be sure to note the discontinuous $x$-axis.

\begin{figure}[H]
    \centering
    \includegraphics[width=0.85\linewidth]{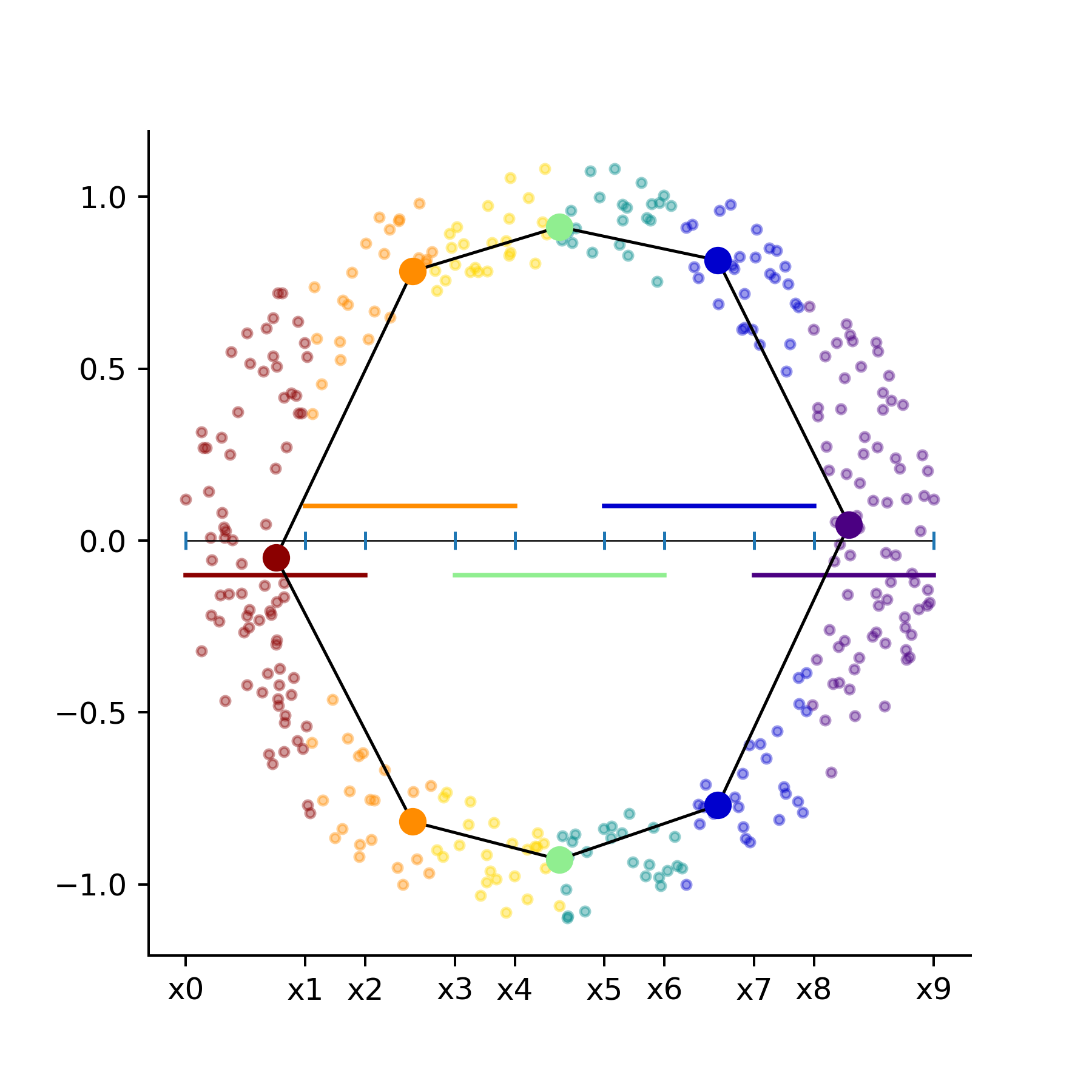}
    \caption{\textit{Simplicial Complex from } $\mathcal{A}'$}
    \label{fig:simplicial}
\end{figure}

Figure \ref{fig:simplicial} is a little busy, so let's break it down. The original data of $\mathcal{A}'$ is there for reference and thus is faded. The $x$-axis is relabeled with the intervals, and the intervals are depicted along the diameter of the annulus parallel to the $x$-axis. The large singular points in each interval are merely the average values of the points in each of these intervals, or the center of mass if each point is weighted equally. The size of the point could vary by the count of data points in the interval it is representing, but here we forego this method. This lower-dimensional approximation of the loop is what becomes particularly useful computationally.

In our example, we had to separate the positive and negative values or else we would lose some information, so upper portions of the circle stay upper, and likewise for the lower. This concept would correspond perhaps with binning the data \textit{along} a circle (and not along the x-axis). Otherwise the positive and negative values would average to a nearly zero $y$ value and we would be compressing the data, losing the loop, which is not desired.

Keep in mind that the filter and intervals here were chosen quite specifically for their graphical properties to demonstrate the concepts, but are usually necessitated by the question driving the investigation, as you will see in the next section where we return to the original data $\mathcal{A}$. Before this however let us introduce some rigour.

Define a continuous function $f:\mathbb{X} \to \R^d$ with $d\geq1$, and as before let $\mathcal{U} = (U_i)_{i\in \mathcal{I}}$. Well, we can take for example that $\mathcal{U}$ is the colored open covering of intervals in figure \ref{fig:simplicial} (the colored line segments). If we pullback $\mathcal{U}$ by $f^{-1}$, we get the overlapping sections of the annulus as seen in figure \ref{fig:overlaps}. These disconnected components are what we average (treating positive and negative clusters separately), and then we obtain the simplicial complex as depicted by the larger vertices and edges in \ref{fig:simplicial}. This is called the \textbf{refined pull back} \cite{chazal}, it is the collection of connected components of the open sets $f^{-1}(U_i)$, for  $i \in \mathcal{I}$.

What we've conducted in this paper is basically the Mapper Algorithm, which is as follows.

\subsection{The Mapper Algorithm}

\noindent \underline{\textbf{Mapper Algorithm}}

\noindent\textbf{Input:} A data set $\mathbb{X}$ with a metric, a function $f:\mathbb{X} \to \R^d$ and a cover $\mathcal{U}$ of $f(\mathbb{X})$.

\noindent\textbf{Algorithm:} For each $U \in \mathcal{U}$, we need to decompose $f^{-1}(U)$ into clusters $C_{1,U},\ldots,C_{U,k_U}$. Next, compute the nerve of the cover of $X$ defined by the $C_{U,i}\in\mathcal{U}$. 

\noindent\textbf{Output:} The simplicial complex preserving topological invariants of the underlying data set $\mathbb{X}$, consisting of a vertex $v_{U,i}$ per cluster $C_{U,i}$, and edges between $v_{U,i}$ and $v_{U',j}$ if and only if $C_{U,i}\cap C_{U',j}\not=\emptyset$. \cite{chazal}

\subsection{Clusters}
In figure \ref{fig:simplicial} we see some sense of symmetry being preserved. While this example is somewhat arbitrary, the preservation of underlying structure in the data is exactly the utility of TDA. This preservation allows us to reduce our analysis from billions or millions of data points (computationally expensive) down to thousands maybe only hundreds, perform the analyses proper, then infer back up to the original clusters whatever our findings. 

\begin{figure}[H]
    \centering
    \includegraphics[width=0.95\linewidth]{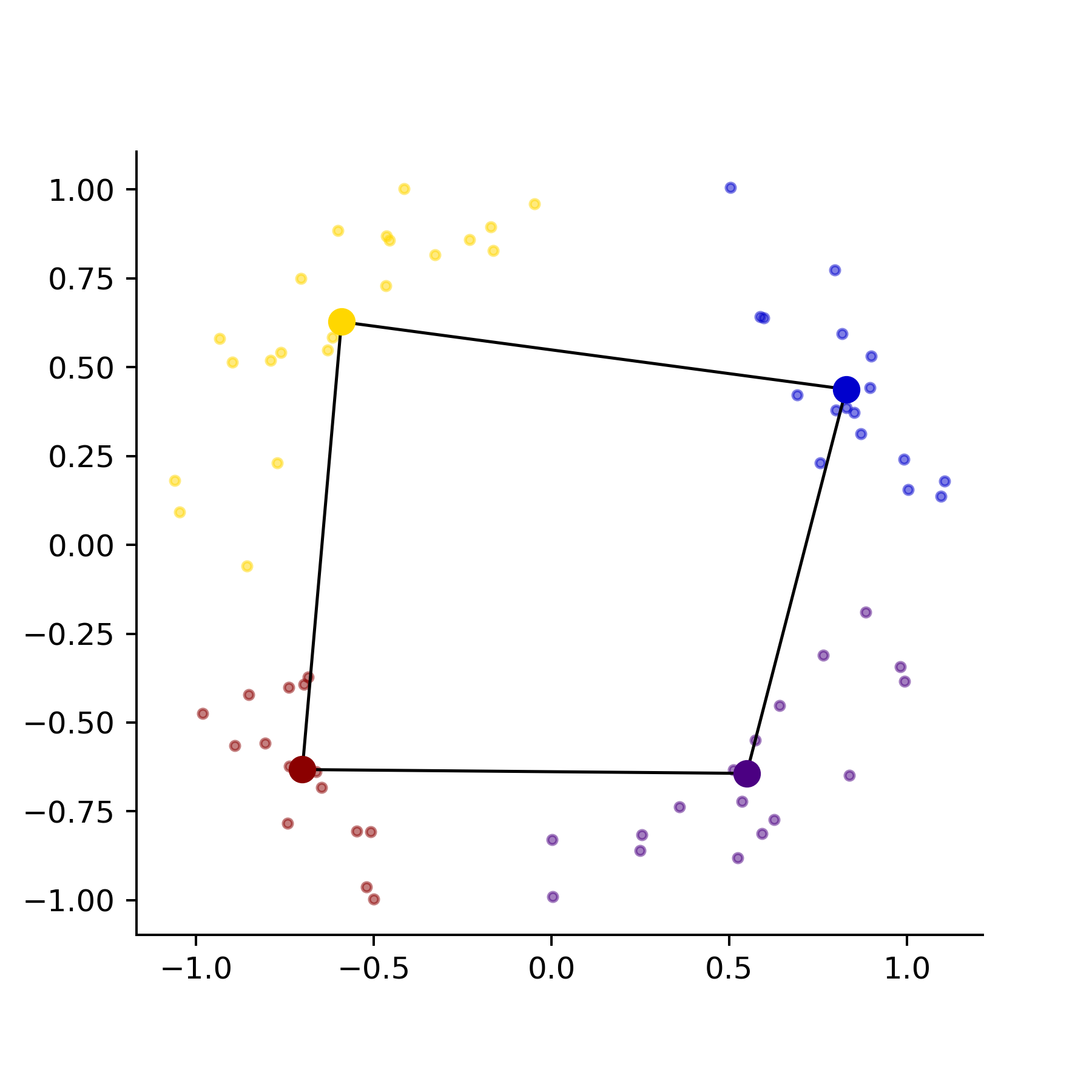}
    \caption{\textit{Simplicial Complex from }$\mathcal{A}$}
    \label{fig:simplicial2}
\end{figure}

Before we partitioned figure \ref{fig:filter}, we introduced the denser set of points $\mathcal{A'}$ for the sake of illustration. Returning to $\mathcal{A}$ we find we can bin (color) by the distinct clusters we found in figure \ref{fig:rad2}. While not a circle per se, our simplex here is still a loop, a `low-resolution' circle as depicted in figure \ref{fig:simplicial2}. This suggests the loop is a topological invariant of the data. In this case we do not need to create a complicated filter as we did for $\mathcal{A}'$, instead choosing to color simply by cluster, then average each cluster, and connecting by separations instead of overlapping intervals. This is not an instance of the mapper algorithm, and is merely a different method of visualizing the relations in our data.

This topological invariance, the preservation of the loops in these examples, can be described by Betti numbers \cite{munch}, thanks to the Nerve Theorem.

\subsection{Betti Numbers and Persistent Homology} \label{betti_numbers}
Since we've established through the Nerve Theorem that a simplicial complex may be homotopy equivalent to a data set $\mathbb{X}$, we need a way to keep track of these geometric shapes. We do this through Betti numbers, which distinguish topological spaces by the connectivity of their $k$-dimensional simplicial complexes. 

\begin{figure}[H]
    \centering
    \includegraphics[width=0.95\linewidth]{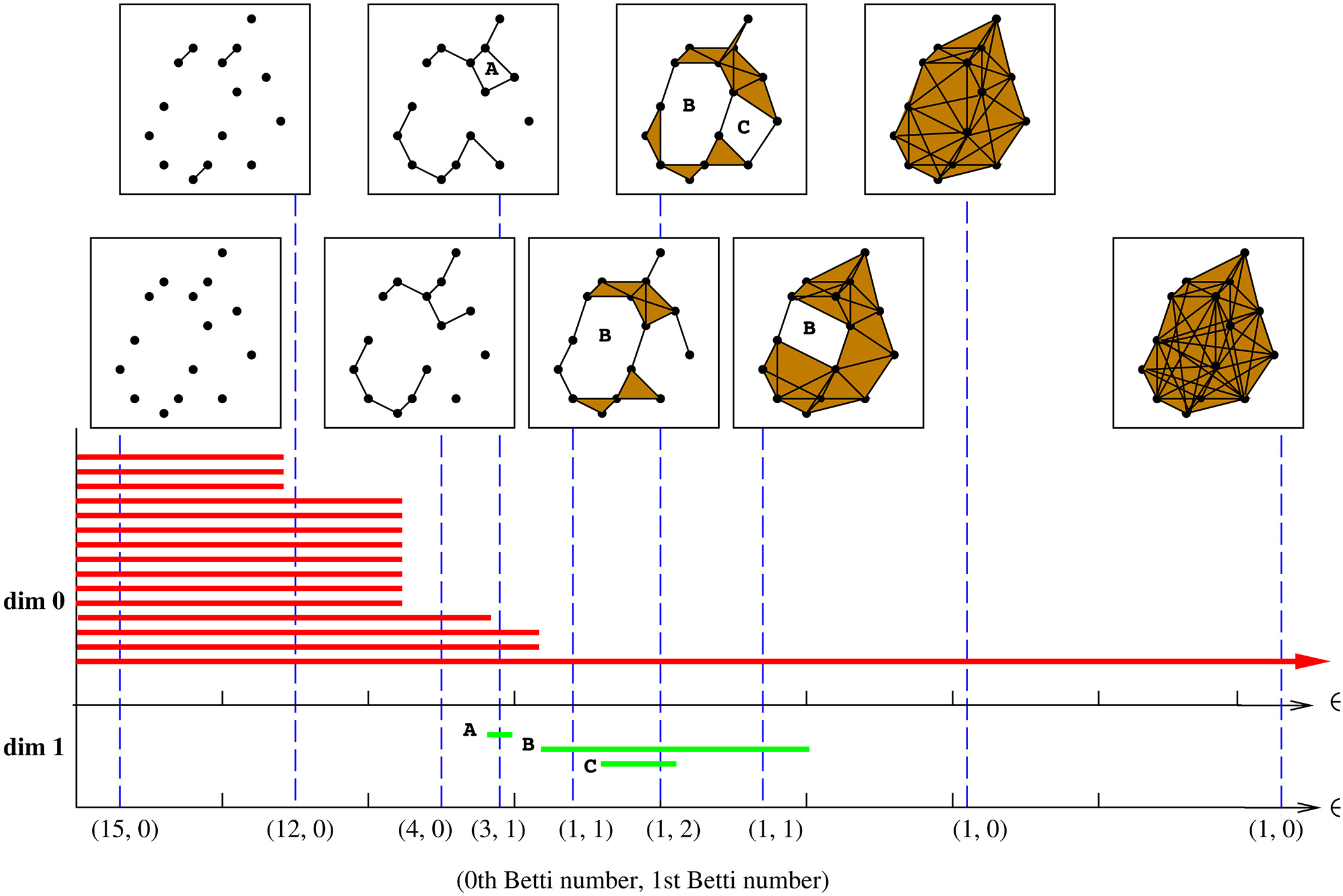}
    \caption{\textit{Persistence Diagram for increasing radii of a Rips complex} - figure via Hylton \textit{et al.} \cite{hylton}}
    \label{fig:persistence}
\end{figure}

First, you may have noticed, there are a lot of choices we've had to make. There is cover, the filter function, the parameter on the homotopy $h(x,\rho)$, et cetera. Indeed, in section \ref{connectivity}, as we increased the radius $r$, loops and holes appeared and disappeared. If we look at all possible values of $r$, we can obtain something called the persistence diagram (see figure \ref{fig:persistence}). 

We saw that between figure \ref{fig:rad3} and \ref{fig:rad4} that the loop structure finally appeared, then by figure \ref{fig:rad6} the loop has disappeared as $r$ approaches $L_\infty$. This appearing and disappearing is the range of persistence for the particular structure of the data, how long (in terms of the radius of the balls) does a structure persist? We're not going to broach this subject too rigorously but can recommend \cite{bubenik} for further reading. It is important to know about this topic in relation to TDA.

To tie together Betti numbers and persistent homology, we note that a \textbf{homology} is divided into dimensions that represent the dimension of the structure \cite{munch}. 0-dimensional homology measures clusters, 1-dimensional homology measures loops, 2-dimensional homology measures holes or voids, and the dimensions can climb beyond this but again we defer to \cite{bubenik}. 

\textbf{Betti numbers} $\beta_k$ for some $k$ are simply the rank of the $k$-dimensional homology group. Below (table \ref{tab:betti}) are listed some simple and common structures with their Betti numbers. Notice a solid disc is homotopy equivalent to a point, and hence their Betti number sequences are identical, they share a homology group, and this is pretty suggestive of the value of this sort of analysis \cite{carlsson}. What you can say about the homotopy equivalence of a data set, you can say about the data set! Remember, millions or billions point data sets are common these days, being able to describe them very well with a couple dozen or a hundred vertices is very, very powerful.

\begin{table}[H]
    \centering
    \begin{tabular}{|c|cccc|} \hline
        \textit{surface}& $\beta_0$ & $\beta_1$ & $\beta_2$ & $\beta_3$\\ \hline
        point           & 1 & 0 & 0 & $0$\\
        solid disc      & 1 & 0 & 0 & $0$\\
        loop            & 1 & 1 & 0 & $0$\\
        sphere          & 1 & 0 & 1 & $0$\\
        torus           & 1 & 2 & 1 & $0$\\
        klein-bottle    & 1 & 1 & 0 & $0$ \\\hline
    \end{tabular}
    \caption{\textit{Betti numbers } $\beta_k$ \textit{ for various surfaces}}
    \label{tab:betti}
\end{table}

And so we examine finally a persistence diagram (figure \ref{fig:persistence} for a data set $\mathbb{X}$, over all radii $0 < r \leq L_\infty$). The horizontal axis represents the chosen radius, and the vertical axis represents the appearance and persistence until disappearance of structures indicated by Betti numbers, of which only $\beta_0$ and $\beta_1$ need be considered for the data $\mathbb{X}$ visualized above the persistence diagram in the same figure.

Notice when loop $A$ arrives in the simplicial complex, there appears in the persistence diagram a short green line where $A$ appears then $disappears$ between the corresponding values of $r$, and it is categorized in $\beta_1$, which characterizes loops! Similarly for loops $B$ and $C$. No voids $\beta_2$ show up is this data set.

There are 15 data points, you'll see 3 short red bars of equal lengths corresponding the the first 3 clusters seen in the simplicial complex. Provided is the pair (12,0) to indicate the three pairs of points ("clusters") and the remaining 9 singletons. As the radius increases you see $A$ appear corresponding to (3,1), where there are three components (two separate connected components and one singleton), along with one cycle/loop. Furthermore, you can see that loop $B$ persists while loop $C$ comes and goes, before $B$ also collapses.

It is these persistence diagrams that guide what radius to use to obtain a desired amount of clusters. It is what reveals to you the structure of your data and how it changes as you investigate its connectivity in it's nearest-neighbors-graph. Further reading is advised to discover more about this wonderful new set of tools.

\section{Discussion}
Topological data analysis provides a new set of tools and methods to get a grasp on your data. It preserves the structure of your data, as seen with the topological invariants of Betti numbers and persistence homology. It has shown to be effective for analysis of complex data sets, from breast cancer to basketball positions and politics, it provides novel insight into questions you may want to address regarding your data. 
Further reading into persistence homology and persistence diagrams will yield you increasingly sophisticated tools for analyzing your data topologically.

If you would like to see the code for the genesis of the custom graphs for this paper they can be viewed at bit.ly/TDA\_2020. The code is customized for the presentation of the random data and likely not transferable to new data sets easily. If you'd like to use the Mapper algorithm there is much literature out there (start with the original \cite{singh}), and there are packages for at least Python and R.

\end{multicols}
\bibliography{mybib}{}
\bibliographystyle{plain}
\end{document}